\documentclass[11pt]{article}
\usepackage{authblk}

\usepackage{geometry}
\usepackage{amsmath}
\usepackage{amssymb}
\usepackage{dsfont}
\usepackage{bbm}
\usepackage{graphicx}
\usepackage{float}
\usepackage{algorithm}
\usepackage{algpseudocode}
\usepackage{listings,lstautogobble}
\usepackage{tikz}
\usepackage{courier}
\usetikzlibrary{decorations.markings}
\usepackage{graphicx}
\usepackage{color}
\usepackage{hyperref}
\hypersetup{
    colorlinks=false, 
    citecolor=black,
     linkcolor=black,
     urlcolor=pink,
    linktoc=all,     
    linkcolor=black,  
}

\numberwithin{equation}{section} 
\usepackage{authblk}

\newcommand\abs[1]{\left|#1\right|}
\newcommand\norm[1]{\left\lVert#1\right\rVert}
\newcommand{\mnorm}[1]{{\left\vert\kern-0.25ex\left\vert\kern-0.25ex\left\vert#1\right\vert\kern-0.25ex\right\vert\kern-0.25ex\right\vert}}
\newcommand\ip[1]{\left\langle#1\right\rangle}
\newcommand{\tr}[1]{\text{tr}\left(#1\right)}
\newcommand{\conv}[1]{\textbf{conv}\left(#1\right)}

\newcommand{\rank}[1]{\text{rank}\left(#1\right)}
\newcommand{\Span}[1]{\textbf{span}\left(#1\right)}

\newcommand{\R}{\mathbb{R}}

\newcommand{\A}{\mathcal{A}}

\newcommand{\ones}{\mathbbm{1}}

\newcommand{\paren}[1]{\left(#1\right)}
\newcommand{\curly}[1]{\left\{#1\right\}}
\newcommand{\brac}[1]{\left[#1\right]}

\newcommand{\supp}[1]{\text{supp}\left(#1\right)}
\newcommand{\qedhere}{\hfill\blacksquare}
\newcommand{\proj}{\text{proj}}

\newcommand{\defeq}{\mathrel{:\mkern-0.25mu=}}

\DeclareMathOperator*{\argmax}{arg\,max}

\title{Sparse Optimization on General Atomic Sets: Greedy and Forward-Backward Algorithms}
\author[1]{Thomas Zhang}
\affil[1]{Yale University \\ \texttt{thomas.zhang@yale.edu}}
\date{}

\begin{document}
\maketitle


\newtheorem{theorem}{Theorem}[section]
\newtheorem{lemma}[theorem]{Lemma}
\newtheorem{proposition}[theorem]{Proposition}
\newtheorem{corollary}[theorem]{Corollary}
\newtheorem{observation}[theorem]{Observation}
\newtheorem{definition}[theorem]{Definition}
\newtheorem{remark}[theorem]{Remark}
\newtheorem{claim}[theorem]{Claim}
\newtheorem{condition}[theorem]{Condition}

\begin{abstract}
    We consider the problem of sparse atomic optimization, where the notion of ``sparsity'' is generalized to meaning some linear combination of few atoms. The definition of atomic set is very broad; popular examples include the standard basis, low-rank matrices, overcomplete dictionaries, permutation matrices, orthogonal matrices, etc. The model of sparse atomic optimization therefore includes problems coming from many fields, including statistics, signal processing, machine learning, computer vision and so on. Specifically, we consider the problem of maximizing a restricted strongly convex (or concave), smooth function restricted to a sparse linear combination of atoms. We extend recent work that establish linear convergence rates of greedy algorithms on restricted strongly concave, smooth functions on sparse vectors to the realm of general atomic sets, where the convergence rate involves a novel quantity: the ``sparse atomic condition number''. This leads to the strongest known multiplicative approximation guarantees for various flavors of greedy algorithms for sparse atomic optimization; in particular, we show that in many settings of interest the greedy algorithm can attain strong approximation guarantees while maintaining sparsity. Furthermore, we introduce a scheme for forward-backward algorithms that achieves the same approximation guarantees. Secondly, we define an alternate notion of weak submodularity, which we show is tightly related to the more familiar version that has been used to prove earlier linear convergence rates. We prove analogous multiplicative approximation guarantees using this alternate weak submodularity, and establish its distinct identity and applications.
    
\end{abstract}


\section{Background and Definitions}

Sparsity in its many forms is central to a variety of problems across statistics and computer science. In general, these problems usually require the estimation of some model whose dimension is much higher than the number of measurements that can be feasibly made. However, if one has the belief or knowledge that the model is constrained in some way that makes it feasibly estimable by few measurements, then sparse optimization becomes a problem of interest. The notion of sparsity differs from problem to problem: in linear least squares, one seeks sparsity in the support of the coefficient vector; in matrix completion, one seeks sparsity in the spectrum of a matrix; in ranked elections, one seeks sparsity in the number of permutations. Atomic sets are an incorporating model that often elegantly capture these different notions of ``sparsity''.

Given an inner product space $\mathcal{H}$, e.g.\ $\R^n$, $M^{m\times n}(\R)$, an atomic set is a (possibly uncountable) set of vectors $\A = \curly{v_i} \subseteq \mathcal{H}$ that is symmetric: if $v \in \A$ then $-v \in \A$. For convenience of notation, many of the atomic sets we mention later are not immediately symmetric; in these cases it is sufficient to assume we are instead dealing with $\A \cup -\A$. We note that the convex hull $\conv{\A}$ contains $0$ and is a polytope when $\A$ is finite. The symmetricity of $\A$ is important for the following property: $\conv{\A}$ induces a norm from its gauge function that we call the ``atomic norm'' induced by $\A$:
\[
\norm{x}_{\A} \defeq \inf \curly{t > 0 : x \in t\cdot \conv{\A}}.
\]
The dual atomic norm is then defined
\[
\norm{x}_{\A^\ast} \defeq \sup_{\norm{z}_{\A} = 1} \ip{z, x}.
\]
Familiar examples of atomic sets include the aforementioned coordinate basis vectors $\curly{e_i} \subset \R^n$ and rank-one matrices $\curly{uv^\top} \subset \R^{n \times d}$. These examples provide a nice intuition to why certain atomic sets induce sparsity: in the coordinate basis vector case, the atomic unit ball is precisely the $\ell^1$ unit ball, which is a polytope, yielding vertex solutions--corresponding to individual atoms--when maximizing/minimizing convex/concave functions. This motivates the study of algorithms for atomic norm regularization \cite{RSW15, CRPW12, LJ15, LKTJ17}. However, in this paper we are concerned with greedy algorithms that \textit{construct} sparse atomic solutions. Such an approach is desirable for multiple reasons. Firstly, atomic norm regularization requires solving a convex program at each iteration; for many atomic sets, atomic norm regularization requires semi-definite programming, which is prohibitively costly for any somewhat high-dimensional problem. Secondly, atomic norm regularization only implicitly induces sparsity, and requires fine-tuning of parameters to get the right degree of sparsity. Thirdly, sharp sufficient conditions under which atomic norm regularization will recover the true sparse solution are in general unknown, and the conditions that are known (e.g.\ Restricted Isometry Property) are often computationally infeasible \cite{TP13}. Therefore, one may be motivated to consider the greedy approach to sparse atomic optimization, where at each iteration the locally optimal (by some metric) atom is added to the active set, thus giving us explicit control over the sparsity of the resulting solution. However, whereas we are guaranteed to converge to the optimal solution in atomic norm regularization, we must face the possibility of a suboptimal solution. One of the main goals of this paper is therefore to establish that the possibly suboptimal greedy solution is comparable to the optimal sparse solution. The problem we consider in this paper is the following ``sparse atomic'' maximization:
\begin{align*}
    (P)\quad \max\,\,&f\paren{\sum_{i=1}^r c_i v_i} \\
    \text{s.t.}\,\,&c_i \in \R, \,\,v_i \in \A, \,\,i = 1,\dots, r.
\end{align*}
We now assume that $f$ is restricted strongly concave and restricted smooth, which are defined as follows.

\begin{definition}[Restricted Strong Concavity, Restricted Smoothness \cite{NRWY12, LW13}] \label{Def:RSCRS}
A function $f: \mathcal{H} 
\to \R$ is restricted strongly concave with parameter $\mu_\Omega$ and restricted smooth with parameter $L_\Omega$ if for all $x,y \in \Omega \subset \mathcal{H}$,
\begin{align*}
    -\frac{\mu_\Omega}{2} \norm{y - x}^2 \geq f(y) - f(x) - \ip{\nabla f(x), y - x} \geq -\frac{L_\Omega}{2} \norm{y - x}^2.
\end{align*}
\end{definition}

We remark that if $\Omega' \subseteq \Omega$, then by first principles
\begin{equation} \label{eq:rscrm}
L_{\Omega'} \leq L_\Omega, \quad \mu_{\Omega'} \geq \mu_{\Omega}.
\end{equation}

In our paper, we often shorten notation and treat $\Omega$ as the whole ambient space, such that restricted strong concavity and restricted smoothness just become their unrestricted counterparts. However, strictly speaking, setting $r$ as the sparsity constraint of $(P)$, then it is sufficient to treat $\Omega$ as the set of all elements of the ambient space that can be written as a linear combination of no more than $2r$ atoms. As a shorthand, we may write the corresponding strong concavity and smoothness parameters as $\mu_{2r}$ and $L_{2r}$. We note that the additional flexibility of ``restricted'' strong concavity and smoothness turns out to be crucial in admitting important problems and objective functions into the model \cite{NRWY12, LW13, NW12, EKDN18}.

Khanna et al.\ \cite{KEDN17, KEDN17sparse, EKDN18} have shown that greedy algorithms attain multiplicative approximations of the optimal solution within $r$ iterations for restricted strongly convex, restricted smooth functions over sparse vectors and low-rank matrices. Our first aim is to show that these algorithms and approximation ratios can be extended in some way to general atomic sets. However, at face value, the definition of ``atomic set'' is extremely broad, and therefore one can easily construct poorly behaving atomic sets where greedy algorithms will not achieve any sort of approximation to the optimal sparse solution in many, many iterations. Therefore, we must introduce a way to measure the structure of an atomic set, in particular its suitability for greedy algorithms. This will come in the form of the ``atomic set condition number'' that will be introduced in the next section. We will end up showing that the approximation guarantee of the greedy algorithm has a very intuitive dependency on three components:
\begin{itemize}
    \item conditioning of the objective function;
    \item conditioning of the underlying atomic set;
    \item the number of greedy steps.
\end{itemize}
The precise meaning of the first two items will be formalized later. This echoes the form of earlier approximation guarantees \cite{LJ15}. However, the concrete notion of ``atomic condition number'' that we introduce has many immediate benefits. Firstly, it generally leads to tighter approximation guarantees than other structural measures of an atomic set. Secondly, the sub-problem of explicitly computing or bounding the atomic condition number is relatively straightforward, as it primarily involves elementary linear algebra, and does not require computing complicated geometric values involving the width or volume of convex bodies. Therefore, if one were to formulate a structured optimization problem in the language of sparse atomic optimization, one could obtain explicit \textit{post hoc} numerical approximation guarantees by deriving the atomic condition numbers. We demonstrate this by deriving a number of atomic condition numbers for common atomic sets in the appendix.

After establishing the importance of atomic condition numbers, our second aim is to revisit and similarly extend to the general atomic setting a notion that has recently been the power tool in establishing greedy approximation guarantees: weak submodularity. One may be familiar with the $\paren{1 + e^{-1}}$ approximation guarantee of the greedy algorithm on non-negative submodular functions \cite{NWF78}. It has been shown in recent work \cite{EKDN18, DK11} that certain families of restricted strongly convex, smooth functions can be transformed into ``weakly'' submodular functions, for which the greedy algorithm attains good approximation guarantees that decay gracefully depending on how ``weakly'' submodular the function is. We show that greedy algorithms also attain nice approximation guarantees in the language of weak submodularity that are distinct from the ones derived using atomic condition numbers: weak submodularity is a notion that has an identity distinct from continuous optimization. We give a motivating example that shows weak submodularity has utility outside the realm of sparse atomic optimization. Roughly speaking, whereas greedy approximation guarantees for sparse atomic optimization have a separate dependence on the conditioning of the atomic set and the conditioning of the objective function, greedy approximation guarantees for weak submodular maximization technically depends only on the weakly submodular function. In previous bounds relating sparse optimization and weak submodularity, one could say the stars aligned and allowed the good conditioning of the objective function and atomic set to translate to a useful notion of weak submodularity.

\section{Atomic Set Properties and Greedy Improvements}

Atomic sets in full generality can be succinctly characterized as any set $\A$, possibly uncountable, in a Hilbert space that is symmetric: if $v \in \A$, then $-v \in \A$. This is the definition used in recent literature concerning the convergence properties of Frank-Wolfe-type algorithms \cite{RSW15, LJ15, GIZ17}: in the past, the term ``atom set'' has been used to refer to certain particular examples of the above definition, in particular the finite-dimensional elementary basis vectors $\curly{e_i} \subset \R^n$ and elements of a dictionary \cite{CDS01, MZ93, DE03, Tr04, GV06}. Examples of atomic sets include the aforementioned instances in addition to rank-one matrices, Dirac measures, orthogonal matrices, permutation matrices, as well as group-sparse atoms \cite{CRPW12, RSW15}.

A unifying goal in introducing the notion of atoms is sparsity. Consider the convex hull $\conv{\A}$: the maximum of a convex function (in particular linear) over it is attained at an extreme point, i.e.\ an atom. This is part of the intuition behind Frank-Wolfe-type algorithms, and also underlies the motivation for atomic regularization \cite{CRPW12, RSW15, LJ15}. In atomic regularization, an objective function of the form $f(x) + \lambda \norm{x}_{\A}$, where $\norm{x}_{\A}$ is the gauge norm induced by $\conv{\A}$, is considered. One may be familiar with the analysis of particular examples of atomic regularization, such as LASSO, where $\A = \curly{e_i}^n$, $\norm{\cdot}_{\A} = \norm{\cdot}_1$ \cite{Ti96, CDS01, CRT06}, and nuclear-norm minimization, where $\A = \curly{uv^\top}$, $\norm{\cdot}_{\A} = \norm{\cdot}_{\ast}$ \cite{RFP10, FHB01}. If $f(x)$ is convex, then gradient-based methods have provably good convergence properties and enjoy Nesterov-type acceleration \cite{N13-p, BT09}, since $f(x) + \lambda \norm{x}_{\A}$ is a convex composite objective. Regardless, there are two immediate drawbacks to the regularization approach to sparse atomic optimization. Firstly, blackbox convex programming solvers are difficult to scale to many important modern problems. For example, matrix completion (low-rank optimization) turns into semidefinite programming, which struggles in practical efficiency as the dimension of the matrix inflates. Secondly, the sparsity induced by the regularization term is implicit: many of the sharper conditions to guarantee a sparse solution are provably expensive to verify \cite{TP13, BDMS13}. This serves as the motivation for a greedy algorithm approach to sparse atomic optimization. In particular, we consider the approach of adding one atom greedily to an active set per iteration, such that we have explicit control over the sparsity. Certainly, the greedy approach may not find the optimal sparse solution, but the ultimate aim is to bound the improvement per iteration to guarantee that the greedy algorithm can generate a sparse solution that is a good approximation of the optimal.

Atomic sets are not created equal. For such a broad definition, we cannot hope to find approximation guarantees that are both strong and universally applicable for sparse atomic optimization. To capture the effectiveness of greedy steps on a given atomic set, let us define the following quantities.

\begin{definition}[Atomic Condition Number]
Let $\A \subseteq \mathcal{H}$ be a given atomic set in inner product space $\mathcal{H}$. We define the atomic condition number of $\A$ to be the largest value such that for any vector $v \in \mathcal{H}$, $v \neq 0$, there exists an atom $a \in \A$ such that
\[
\frac{\abs{\ip{v,a}}}{\norm{v} \norm{a}} \geq \theta.
\]
In other words,
\[
\theta \defeq \min_{\norm{v} = 1} \max_{a \in \A} \frac{\abs{\ip{v,a}}}{\norm{a}}.
\]
\end{definition}

\begin{definition}[Sparse Atomic Condition Number]
Let $\A \subseteq \mathcal{H}$ be a given atomic set in inner product space $\mathcal{H}$. Given $L \subset \A$, let the atomic condition number with respect to $L$, $\theta(L) \geq 0$, to be
\[
\theta(L) \defeq \min_{\substack{\norm{v} = 1 \\ v \in \Span{L}}} \max_{a \in \A} \frac{\abs{\ip{v,a}}}{\norm{a}}.
\]
The $r$-sparse, or simply sparse, atomic condition number is defined as the minimum atomic condition number over all $r$-subsets: $\theta_r \defeq \min_{\abs{L} \leq r} \theta(L)$. In other words,
\[
\theta_r \defeq \min_{\substack{L \subseteq \A \\ \abs{L} \leq r}} \min_{\substack{\norm{v} = 1 \\ v \in \Span{L}}} \max_{a \in \A} \frac{\abs{\ip{v,a}}}{\norm{a}},
\]
\end{definition}
where $\abs{L}$ indicates the cardinality of $L$.

Essentially, the atomic condition number measures how dense an atomic set is. A large $\theta$ means every vector is reasonably close to an atom. We can identify some representative atomic sets with properties that lead to tight lower bounds on $\theta_r$. We will see that larger values of $\theta_r$ directly correspond to stronger bounds on greedy improvement.

\begin{itemize}
    \item (Topologically) dense atomic sets, e.g.\ Euclidean unit sphere. In that case, $\theta = 1$ trivially, since any vector is well-approximated by an atom. That said, sparsity with respect to these atomic sets is not often very meaningful.

    \item Orthogonal basis. $\theta_r \geq r^{-1/2}$. This class of atomic sets is perhaps the nicest ``meaningful'' atomic set, particularly due to the following property: given two sets of atoms $S$ and $T$, where $\abs{T} = k$, if we define $T'$ as the minimal set of atoms such that
    \[
    \Span{T'} \supseteq \proj_{S^{\perp}}(T),
    \]
    we observe $\abs{T'} \leq k$, as $T'$ is simply $T$ minus the atoms it shares with $S$. This is a key property used in Elenberg et al.\ \cite{EKDN18} to prove strong approximation guarantees for greedy algorithms where $\A$ is the set of elementary basis vectors (sparse optimization). This is a property lost even when considering union of orthogonal bases (dictionary learning).

    \item Other atomic sets with $\theta_r$ dependent only on $r$, for example rank-one matrices: $\theta_r \geq r^{-1/2}$ as the linear combination of $r$ rank-one matrices is at most rank $r$. We observe that rank-one matrices do not satisfy the additional property mentioned above.

    \item Atomic sets with $\theta_r$ bounded away from $0$, but possibly dependent on the ambient dimension. For example, the set of orthogonal matrices $U \in M_n$, where $\theta_r = n^{-1/2}$. The proof of this bound will be relegated to the appendix.

\end{itemize}

Table~\ref{sparse theta values}
contains a larger sample of atomic sets and their respective lower bounds on the atomic condition numbers. The proofs of the bounds in Table~\ref{sparse theta values} can be found in the appendix.

\begin{table}[] 
\begin{tabular}{|l|l|l|}
\hline
\textbf{Atomic set}  & \textbf{Atomic cond.\ number} $\theta$              & \textbf{Sparse} $\theta_r$ \textbf{value} \\ \hline
Standard basis $\curly{e_i}_{i=1}^n \subset \R^n$           & $n^{-1/2}$                & $r^{-1/2}$     \\ \hline
$m \times n$ rank-one matrices $uv^\top$   & $(\min\curly{m,n})^{-1/2}$  & $r^{-1/2}$     \\ \hline
Disjoint group-sparse atoms $\mathcal{P}\paren{\curly{e_i}_{i=1}^n}$      & $(\text{\# groups})^{-1/2}$ & $r^{-1/2}$     \\ \hline
2-ortho basis $\curly{\phi_i}^n \cup \curly{\psi_j}^n \subset \R^n$ &  $\Omega\paren{n^{-1/2}}$         & $\Omega\paren{r^{-1/2}}$ when $r \leq \mu\paren{\A}^{-1}$        \\ \hline
Binary sign vectors $\curly{\pm 1}^n$   & $n^{-1/2}$ &  $n^{-1/2}$  \\ \hline
$n \times n$ orthogonal matrices  & $n^{-1/2}$                & $n^{-1/2}$     \\ \hline

\end{tabular}
	\caption{Sparse Atomic Condition Numbers} \label{sparse theta values}
\end{table}

We note that even though adding more vectors to an atomic set can only increase $\theta$, and $\theta_r \geq \theta$, it is not necessarily the case that $\theta_r$ monotonically increases with the addition of vectors to the atomic set. In fact, the gap between $\theta_r$ and $\theta$, if it exists, can be made arbitrarily small by adding just one vector to the atomic set. Consider the following simple example. Let $\A = \curly{e_i}_{i=1}^n$. Consider adding to $\A$ the vector $v = e_1 + \varepsilon \ones$, where $\ones$ is the all-ones vector, and $\varepsilon$ is an arbitrary small value. Then for any $r \geq 2$, we consider the subset containing $e_1$ and $v$. The span of that subset will contain $n^{-1}\ones$, which satisfies
\[
\max_{a \in \A \cup \curly{v}} \frac{\abs{\ip{\ones, a}}}{\sqrt{n} \norm{a}} \leq \frac{1}{\sqrt{n}} + \delta,
\]
where $\delta$ can be made arbitrarily small by shrinking $\varepsilon$. Therefore, we have shown that we can make the sparse atomic condition number $\theta_r$ for the standard basis, which is normally $r^{-1/2}$, arbitrarily close to its atomic condition number $\theta = n^{-1/2}$ simply by corrupting the atomic set with one vector. We will later see that this phenomenon is closely tied to how well a greedy approach can find good sparse solutions. In short, the structure of an atomic set is important!

\subsection{Atomic Sets and Greedy Optimization}

The atomic condition numbers have a direct relationship with bounds on greedy algorithmic performance. First, we state the following combinatorial reformulation of problem $(P)$ on which we apply greedy algorithms:
\begin{align*}
(P) \quad \max &\,\, g(U) \\
\text{s.t.} &\,\, U \subset \A,\,\,\abs{U} \leq r.
\end{align*}
In the cases of our concern, $g: 2^{\A} \to \R$ is a set function defined on subsets of the atomic set:
\[
g(U) \defeq \max_{x \in \Span{U}} f(x) - f(0),
\]
where $f(x)$ is some (restricted) strongly concave and smooth function that we want to maximize. Defining $g(\emptyset) = 0$, we observe that $g$ is a non-negative function. Let $g(U^*) = f(x^*) - f(0)$ be the optimal value of $(P)$. Similarly, given set $U$, we define
\begin{align*}
B^{(U)} \defeq \argmax_{x \in \Span{U}} f(x) - f(0).
\end{align*}
In other words, $B^{(U)}$ is the vector in $\mathcal{H}$ such that $f(B^{(U)}) = g(U) + f(0)$. We have the following standard lemma relating the norm of the gradient at a given point to the objective value.

\begin{lemma}\label{Lem:GradDistOpt}
Let $x^*$ be the optimal solution to $(P)$. Then for any $x$, we have
\begin{align*}
    \frac{\norm{\nabla f(x)}^2}{2L} \leq f(x^*) - f(x) \leq \frac{\norm{\nabla f(x)}^2}{2\mu}.
\end{align*}
\end{lemma}

\textit{Proof:} by the concavity of $f$, we have
\begin{align*}
    f\paren{x^*} &\leq f(x) + \nabla f(x)^\top \paren{x^* - x} - \frac{\mu}{2}\norm{x^* - x}^2 \\
    &\leq f(x) + \norm{\nabla f(x)} \norm{x^* - x} - \frac{\mu}{2}\norm{x^* - x}^2 \\
    &\leq f(x) + \frac{\norm{\nabla f(x)}^2}{2\mu}.
\end{align*}

On the other hand, from the smoothness of $f$, we have
\begin{align*}
    f\paren{x^*} &\geq f\paren{x + \frac{1}{L}\nabla f(x)} \\
    &\geq f(x) + \frac{1}{L}\nabla f(x)^\top \nabla f(x) - \frac{1}{2L}\norm{\nabla f(x)}^2 \\
    &= f(x) + \frac{\norm{\nabla f(x)}^2}{2L}.
\end{align*}
$\qedhere$

Therefore, given any subset $U$, if we pick an atom $v$ satisfying
\[
\frac{\ip{\nabla f(x_U), v}}{\norm{\nabla f(x_U)} \norm{v}} \geq \theta,
\]
where such a $v$ is guaranteed to exist by the definition of $\theta$, we can lower bound the gain from $g(U)$ to $g(U \cup \curly{v})$. From now on, assume that the elements of the atomic set $\A$ are normalized: $\norm{v} = 1$ for all $v \in \A$.
\begin{align*}
    g\paren{U \cup \curly{v}} &\geq \max_t f\paren{x_U^* + tv} - f(0)\\
    &\geq f\paren{x_U^*} + \max_t\brac{t \nabla f\paren{x_U^*}^\top v - t^2\frac{L}{2} \norm{v}^2} - f(0)\\
    &\geq g(U) + \theta^2 \frac{\norm{\nabla f\paren{x_U^*}}^2}{2L}, \\
    g\paren{U \cup \curly{v}} - g\paren{U} &\geq \theta^2 \frac{\norm{\nabla f\paren{x_U^*}}^2}{2L} \\
    &\geq \theta^2 \frac{\mu}{L} \paren{g(U^*) - g(U)}
\end{align*}
by Lemma \ref{Lem:GradDistOpt}. In other words, greedily adding an atom yields an objective gain toward the optimal that can be lower bounded multiplicatively. We now introduce the simple greedy algorithm.

\begin{algorithm}
\caption{Greedy$(\A, f, r,\beta)$}
\begin{algorithmic}[1]
\State $U_0 \leftarrow \emptyset$
\For{$t = 1,\dots,r$}
\State Compute $B^{(U_{t-1})}$, $\nabla f\paren{B^{(U_{t-1})}}$
\State $v_t \leftarrow \texttt{PureGreedy}(\A, f, U_{t-1}, \beta)$
\State $\paren{v_t \leftarrow \texttt{OMPSel}(\A, f, U_{t-1}, \beta)}$
\State $U_t \leftarrow U_{t-1} \cup \curly{v_t}$
\EndFor
\State \Return $U_{r}, B^{(U_r)}, f(U_r)$
\end{algorithmic}
\end{algorithm}

\texttt{PureGreedy}$(\A, f, U_{t-1}, \beta)$ is an oracle that returns an atom $v_t$ such that
\[
g\paren{U_{t-1} \cup \curly{v_t}} - g\paren{U_{t-1}} \geq \beta \max_{v \in \A} \paren{g\paren{U_{t-1} \cup \curly{v}} - g\paren{U_{t-1}}}.
\]
On the other hand, \texttt{OMPSel}$(\A, f, U_{t-1}, \beta)$ is a linear oracle that returns an atom $v_t$ such that
\[
\abs{\ip{\nabla f\paren{B^{(U_{t-1})}}, v_t}} \geq \beta \max_{v \in \A} \abs{\ip{\nabla f\paren{B^{(U_{t-1})}}, v}}.
\]
We note that either oracle can be used without affecting any our approximation guarantees. However, \texttt{OMPSel} is usually the more computationally feasible option.

Applying the greedy algorithm to $(P)$, from our lower bound on the gain of greedy addition, the greedy algorithm attains the following multiplicative approximation of the optimal solution.
\begin{align*}
    g(U^*) - g(U_{t}) &\leq \paren{1 - \beta\theta^2 \frac{\mu}{L}}\paren{g(U^*) - g(U_{t-1})} \\
    &\leq \paren{1 - \beta\theta^2 \frac{\mu}{L}}^{t}\paren{g(U^*) - g(\emptyset)} \\
    &\leq \paren{1 - \beta\theta^2 \frac{\mu}{L}}^{t} g(U^*) \\
    &\leq \exp\paren{-\beta\theta^2 \frac{\mu}{L}} g(U^*) \\
    \implies g(U_{t}) &\geq \paren{1 - \exp\paren{- \beta \theta^2 t \frac{\mu}{L}}} g(U^*).
\end{align*}

Up to this point, analogous convergence rates have been shown for other greedy-type methods \cite{LJ15, LKTJ17}: substituting the $\theta = \frac{1}{\sqrt{n}}$ value for sparse optimization, we get an approximation guarantee of the form
\[
g(L_{t}) \geq \paren{1 - \exp\paren{- \beta \frac{t}{n} \frac{\mu}{L}}} g(U^*).
\]
According to this approximation guarantee, if we have run the greedy algorithm for $O(n)$ iterations, we get an approximation ratio solely dependent on the condition number $\mu/L$ and the precision constant $\beta$. However, this is not the end-goal of sparse optimization, as the solution will have $O(n)$ non-zero entries.

\subsection{Tightening Greedy Bound}
As previewed earlier, our goal is to create a bound on the greedy performance that attains a ``constant-factor'' approximation ratio (that is, solely dependent on the condition number and the precision constant) \textit{while} maintaining sparsity of the solution. Here we will show that we can replace $\theta$ in our earlier bounds with $\theta_{2r}$. Assume we have applied the greedy algorithm on $(P)$ for $r$ iterations and attained atom set $U_r$, and that the optimal solution to $(P)$ is $U^*$. Define $V = U_r \cup U^*$. We consider a restricted version of $(P)$:
\begin{align*}
(P_R) \quad \max &\,\, \tilde{f}(x) \\
\text{s.t.} &\,\, x \in \Span{U} \\
&\,\, U \subseteq V,\,\,\abs{U} \leq r,
\end{align*}

where $\tilde{f}(x) \defeq f\paren{\proj_V(x)}$, where $f$ is the restricted strongly concave, smooth objective function in $(P)$. $\tilde{f}(x)$ is concave, as the projection operator is linear. Additionally, $\tilde{f}(x)$ also inherits the strong concavity and smoothness of $f$, as long as $x$ is restricted to $\Span{V}$. In other words, if $f$ is $m_{\Omega}$-restricted strongly concave and $M_{\Omega}$-smooth on $\Omega$, then $\tilde{f}(x)$ is $m_{\Omega}$-restricted strongly concave and $M_{\Omega}$-smooth on $\Omega \cap \Span{V}$. Let us introduce the sparse condition number of the objective function $f$:
\[
\sigma_r \defeq \min_{\substack{L \subseteq \A \\ \abs{L} \leq r}} \min_{u \in \Span{L}} \frac{\mu(u)}{L(u)}.
\]
In other words, $\sigma_r$ is the condition number of the function $f$ over all subspaces of dimension at most $r$. Observe that $\sigma_r \geq \sigma \defeq \mu/L$. Recalling our notation $\mu_{r}$ and $L_r$, indicating the restricted strong convexity and smoothness constants over all subspaces of dimension at most $r$, we also have $\sigma_r \geq \mu_r / L_r$. We note that the latter expression may be more practical to estimate.

Observe that the optimal values of $(P)$ and $(P_R)$ are the same. However, since the dimension of the search space of $(P_R)$ is at most $2r$, we can replace the $\theta$ in previous derivations with $\theta_{2r}$, and $\mu/L$ with $\sigma_r$. We note that $\theta_{2r}$ may not always have a polynomial dependence on $r$; in some cases, $\theta_{2r}$ might be no better than $\theta$, as one can see from the table of atomic condition numbers. The greedy algorithm applied on $(P_R)$, Greedy$(V, \tilde{f}, r, \beta)$, therefore has a convergence rate of
\[
\tilde{g}(U_{t}) \geq \paren{1 - \exp\paren{- \beta \theta_{2r}^2 t \sigma_{2r}}} \tilde{g}(U^*).
\]
Since $\tilde{g}(U^*) = g(U^*)$, if we show that the iterates of Greedy$(V, \tilde{f}, r, \beta)$ are identical with the iterates of the greedy algorithm applied to $(P)$, Greedy$(\A, f, r, \beta)$, then the above convergence rate is actually the convergence of the greedy algorithm on $(P)$. If \texttt{PureGreedy} is used, this is trivial, since $\tilde{f}$ is a restriction of $f$, and all the locally optimal choices for the greedy algorithm on $(P)$ are available in the search space of $(P_R)$. If \texttt{OMPSel} is used instead, the iterates are still identical by applying the chain rule: let $P_V$ denote the projection matrix projecting onto $\Span{V} = \Span{U_r \cup U^*}$ such that $\tilde{f}(x) = f\paren{P_V x}$. By the chain rule we have
\[
\nabla\tilde{f}(x) = P_V \nabla f\paren{P_V x}
\]
Each $v_t$ chosen by the greedy algorithm on $(P)$ satisfies
\[
\ip{\nabla f\paren{B^{(U_{t-1})}}, v_t} = \max_{v \in \A} \ip{\nabla f\paren{B^{(U_{t-1})}}, v}.
\]
By definition of $U_r$, we have $v_t \in U_r$ for all $t = 1,\dots,r$. Since $B^{(U_{t-1})} \in \Span{V}$, $v_t$ also satisfies
\begin{align*}
\ip{\nabla f\paren{B^{(U_{t-1})}}, v_t} &= \max_{v \in \A} \ip{\nabla f\paren{B^{(U_{t-1})}}, P_V v} \\
&= \max_{v \in \A} \ip{P_V \nabla f\paren{P_V B^{(U_{t-1})}}, v} \\
&= \max_{v \in \A} \ip{\nabla \tilde{f}\paren{B^{(U_{t-1})}}, v}.
\end{align*}
Therefore, the locally optimal atom at each iteration on $(P)$ as decided by \texttt{OMPSel} agrees with the locally optimal atom on $(P_R)$. By induction, this implies that the iterates of the greedy algorithm on $(P)$ agree with the iterates of the greedy algorithm on $(P_R)$. Therefore, we have established the improved approximation guarantee.

\begin{theorem} \label{Thm:GreedyFastApproxRatio}
Greedy$(\A, f, k, \beta)$ has the following multiplicative improvement ratio and approximation guarantee:
\begin{align*}
    g(U^*) - g(L_{t}) &\leq \paren{1 - \beta\theta_{2r}^2 \sigma_{2r}}\paren{g(U^*) - g(L_{t-1})} \\
    g(L_{t}) &\geq \paren{1 - \exp\paren{- \beta \theta_{2r}^2 t \sigma_{2r}}} g(U^*).
\end{align*}
\end{theorem}

Referring to our earlier discussion of $\theta_{r}$ values for particular atomic sets, we have the following examples of improved greedy approximation guarantees.

\begin{corollary}[Greedy Feature Selection Convergence Rate]
Consider problem $(P)$, where $\A = \curly{e_i}_{i=1}^n \subset \R^n$. Given a function $f$ (and corresponding function $g$) that satisfies RSC-RS, we have the following lower bound for the performance of the greedy algorithm
\[
g(U_{t}) \geq \paren{1 - \exp\paren{-\beta \theta_{2r}^2 t \sigma_{2r} }}g(U^*) \geq \paren{1 - \exp\paren{-\beta \frac{t}{2r}\frac{\mu_{2r}}{L_{2r}}}}g(U^*)
\]
\end{corollary}

\begin{corollary}[Greedy Low-Rank Optimization Convergence Rate]
Consider problem $(P)$, where $\A = \curly{uv^\top} \subset \R^{m \times n}$. Given a function $f$ (and corresponding function $g$) that satisfies RSC-RS, we have the following lower bound for the performance of the greedy algorithm
\[
g(U_{t}) \geq \paren{1 - \exp\paren{-\beta \theta_{2r}^2 t \sigma_{2r} }}g(U^*) \geq \paren{1 - \exp\paren{-\beta \frac{t}{2r}\frac{\mu_{2r}}{L_{2r}}}}g(U^*)
\]
\end{corollary}

Note that the above approximation guarantees imply that given the sparsity constraint $r$, the greedy algorithm will find a constant-factor approximation of the optimal $r$-sparse solution within $O(r)$ iterations, instead of $O(n)$ iterations. These agree with the current best approximation guarantees (up to small constant factors) for greedy-type algorithms for the above settings \cite{EKDN18, KEDN17, KEDN17sparse, GIZ17, WLLFDY14}.

\subsection{Forward-Backward Schemes}

Our goal is to extend an analogous approximation guarantee to a flexible family of forward-backward algorithms. The motivation of forward-backward algorithms is that ``bad atoms'' contributing little to the objective chosen earlier by the myopic forward steps may be removed later by backward steps to improve the quality of the sparse solution. At its core, the forward-backward paradigm is heuristic, and thus bounds on its performance even in familiar settings and for popular objectives are few and far between, despite bounds existing on the forward-only procedure. We will show that a large class of forward-backward schemes have approximation guarantees no worse than the corresponding forward-only scheme. We propose the following framework for the Forward-Backward scheme.

\begin{algorithm}
\caption{FoBa$(\A, f, k,\beta, \nu)$}
\begin{algorithmic}[1]
\State // $\nu$ is a thresholding constant $0 \leq \nu < 1$
\State $S_0 \leftarrow \emptyset$
\State $t \leftarrow 0$
\While{$\abs{S_t} < k$}
\State $v_{t+1} \leftarrow \texttt{PureGreedy}(\A, f, S_{t}, \beta)$
\State $\paren{v_{t+1} \leftarrow \texttt{OMPSel}(\A, f, S_{t}, \beta)}$
\State $S_{t+1} \leftarrow S_{t} \cup \curly{v_{t+1}}$
\State $t \leftarrow t+1$
\State // Optional: additional conditions to enter backward step
\State $\texttt{Dmg} \leftarrow 0$ \quad // variable tracking the damage done by backward steps
\State $d^+ \leftarrow g(S_t) - g(S_{t-1})$
\While{\textbf{true}}
\State $p \leftarrow \argmax_{v \in \A} g(S_t \setminus \curly{p})$ \quad // picking the element contributing the least
\State $d^- \leftarrow g(S_t) - g(S_t \setminus \curly{p})$
\If{$(\texttt{Dmg} + d^-) < \nu d^+ $}
\State $S_{t+1} \leftarrow S_{t} \setminus \curly{p}$
\State $t \leftarrow t + 1$
\State $\texttt{Dmg} \leftarrow \texttt{Dmg} + d^-$
\Else
\State \textbf{break}
\EndIf
\EndWhile
\EndWhile
\State \Return $S_t, g(S_t)$
\end{algorithmic}
\end{algorithm}

\begin{theorem}[FoBa Convergence Rate] \label{Thm:FoBaConv}
Let $x_S^*$ denote the optimal solution to $(P)$, with sparsity constraint $k$. If for some $0 < c \leq 1$ the forward-only procedure Greedy$\paren{\A, f, k, \beta}$ satisfies a convergence rate of the form
\[
g(S^*) - g(S_{t+1}) \geq \paren{1 - c}(g(S^*) - g(S_t)),
\]
then FoBa$(\A, f, k,\beta, \nu)$ has a convergence rate
\[
g(S_{t}) \geq \paren{1 - \paren{1 - \frac{c}{k}}^{\abs{S_t}}} g(S^*) \geq \paren{1 - \exp(-c \abs{S_t})}g(S^*)
\]
Note that this bound is independent of the thresholding constant $0 \leq \nu < 1$.
\end{theorem}

\textit{Proof of Theorem \ref{Thm:FoBaConv}}: We use induction. First we verify the base cases: $t = 0 \implies S_0 = \emptyset$. By our algorithm, $t = 1$ must be a forward step, and therefore the bound is true due to our assumption of the forward-only convergence rate. Assume the induction hypothesis: at steps $j < t + 1$ we have
\[
g(S_j) \geq \paren{1 - \paren{1 - \frac{c}{k}}^{\abs{S_j}}} g(S^*).
\]
After step $t$, we are at one of the following two scenarios.
\begin{itemize}
    \item Case 1: step $t + 1$ will be a forward step. Therefore, $\abs{S_{t+1}} = \abs{S_t} + 1$. By the same argument made in the proof of Theorem \ref{thm:wksub2}, we have
    \[
    g(S_{t+1}) \geq \paren{1 - \paren{1 - \frac{c}{k}}^{\abs{S_{t+1}}}} g(S^*)
    \]

    \item Case 2: step $t + 1$ will be a backward step. Since we can only take a backward step after making at least one forward step, say our last forward step was at step $t-i, \,i \geq 0$. By the thresholding in the backward step, we have that
    \[
    g\paren{S_{t+1}} \geq g\paren{S_{t-i-1}},
    \]
    and since all steps since $t-i$ are backward steps, we have
    \[
    \abs{S_{t+1}} = \abs{S_{t-i-1}} - i \leq \abs{S_{t-i-1}}.
    \]
    By the induction hypothesis we have
    \begin{align*}
    g\paren{S_{t+1}} &\geq g\paren{S_{t-i-1}} \\
    &\geq \paren{1 - \paren{1 - \frac{c}{k}}^{\abs{S_{t-i-1}}}} g(S^*) \\
    &\geq \paren{1 - \paren{1 - \frac{c}{k}}^{\abs{S_{t+1}}}} g(S^*).
    \end{align*}

    Therefore, when the algorithm terminates at step $N$, we have $\abs{S_N} = k$ and thus
    \[
    g\paren{S_N} \geq \paren{1 - \paren{1 - \frac{c}{k}}^{k}} g(S^*) \geq \paren{1 - \exp(-c)}g(S^*).
    \]
    $\qedhere$

\end{itemize}

Substituting $c = \beta \frac{\mu}{L} \theta_{2k}^2$, we recover the exact same approximation guarantee for the forward-backward scheme as the purely greedy scheme. We cannot hope for a better guarantee in general, as it is possible for no backward steps to have been taken. To some degree, it is also not surprising that the forward-backward scheme is ``as good'' as the purely greedy scheme, but we note that the qualifications and conditions to enter the backward phase and/or to take a backward step can be modified in numerous ways to get better empirical results and will likely still result in similar approximation guarantees; we are proposing but one popular class of forward-backward schemes \cite{Zh11}.

\section{Weak Submodularity}

Recently, strong multiplicative bounds for greedy performance on particular atomic sets (i.e.\ $\A = \curly{e_i}^n$, $\A = \curly{u v^\top}$) were established using the notion of weak submodularity \cite{EKDN18, KEDN17, KEDN17sparse}. In the previous section, we have recovered and extended these bounds to the general atomic setting independent of weak submodularity. However, we note that weak submodularity in the aforementioned papers served predominantly as a convenient combinatorial interpretation of a continuous convex problem. Namely, a notion known as ``submodularity ratio'' \cite{DK11} is developed, and it is essentially shown that a $\mu$-strongly convex, $L$-smooth function can be turned into a $\mu/L$-weakly submodular set function, for which greedy maximization attains a $\paren{1 - \exp(-\mu/L)}$-approximation guarantee, reminiscent of the $\paren{1 - 1/e}$ guarantee for greedy maximization of submodular functions in the seminal paper by Nemhauser et al.\ \cite{NWF78}. We note that the ground sets of the weakly-submodular functions in the aforementioned literature are limited to the standard basis \cite{KEDN17, KEDN17sparse, DK11} or the rank-one matrices \cite{EKDN18}. We can generalize bounds on weak submodularity to general atomic sets that satisfy a certain conditions.

We note that there are rich classes of set functions that cannot be fruitfully converted into the language of convexity, but weak submodularity may be ascribed. In this section, we develop a new notion of weak submodularity which we will prove is interchangeable with versions introduced in earlier papers. Weak submodularity allows us to establish strong multiplicative bounds for the greedy algorithm, and hence the forward-backward algorithm. We define the \textit{subset} submodularity ratio, which is related to what will be called the \textit{disjoint} submodularity ratio seen in earlier literature \cite{DK11, EKDN18, KEDN17}.

\begin{definition}[Disjoint Submodularity Ratio \cite{DK11}] \label{Def:WeakSubm1}
Let $P,Q \subset [p]$ be two disjoint index sets, and $g : [p] \to \R$. The \emph{disjoint} submodularity ratio of $P$ with respect to $Q$ is defined as
\[
\gamma_{P,Q} := \frac{\sum_{i \in Q} \paren{g\paren{P \cup \curly{i}} - g(P)}}{g(P\cup Q) - g(P)}.
\]
\end{definition}

Intuitively, the disjoint submodularity ratio measures the diminishing marginal returns property. The additional adjective ``disjoint'', which is not seen in earlier literature \cite{DK11, EKDN18}, is introduced here to distinguish it from a separate notion of submodularity ratio we define next. Note that $\gamma_{P,Q} \geq 1$ if $f$ is submodular. We can also define the disjoint submodularity ratio of a set $U \subseteq [p]$ with respect to an integer $k > 0$:
\[
\gamma_{U,k} \defeq \min_{\substack{P,Q: P \cap Q = \emptyset;\\ P \subseteq U \\\abs{Q} \leq k}} \gamma_{P,Q}.
\]

Elenberg et al.\ \cite{EKDN18} recently demonstrated in the standard basis setting that functions for satisfying RSC-RS are weakly submodular, where the disjoint submodularity ratio is essentially lower bounded by the restricted condition number: $\gamma_{U, k} \geq \mu_{\abs{U} + k} / L_{\abs{U}+k}$. This leads to a greedy convergence rate of the form
\[
g(P_t) \geq \paren{1 - \exp\paren{-\frac{\mu_{\abs{U} + k}}{L_{\abs{U}+k}}}} g(P^*),
\]
which is precisely analogous to the bound we established for the general atomic setting in the previous section. The reader is directed to the paper by Elenberg et al.\ \cite{EKDN18} or the paper by Das and Kempe \cite{DK11} for more details. A persistent barrier in simply extending their results to connect general sparse atomic optimization and weak submodularity is summarized by the following: not every atomic set contains an orthogonal basis, let alone an orthogonal basis that can be formed starting with any arbitrary element in the atomic set.

Let us introduce a new notion of weak submodularity below, which directly relies on the sense of approximately diminishing marginal returns for a set function. 

\begin{definition}[Subset Submodularity Ratio] \label{Def:WeakSubm2}
Let $U,V \subset [p]$ be index sets such that $U \subseteq V$, and $g: [p] \to \R$. The \emph{subset} submodularity ratio of $U$ with respect to $V$ is defined as
\[
\kappa_{U,V} \defeq \min_{i \in [p]\setminus V} \frac{g\paren{U \cup \curly{i}} - g\paren{U}}{g\paren{V \cup \curly{i}} - g\paren{V}}
\]
\end{definition}
Similarly, we may define the subset submodularity ratio of a set $U \subseteq [p]$ with respect to an integer $k > 0$:
\[
\kappa_{U, k} \defeq \min_{\substack{T,V : T \subseteq U,\, U \subseteq V \\ \abs{V \setminus U} \leq k}} \kappa_{T,V} . 
\]

Now that we have stated results for the disjoint submodularity ratio, it is interesting to note that the disjoint and subset submodularity ratios are tightly related.

\begin{theorem}\label{thm:type1totype2} Given $\gamma = \gamma_{U, k}$ and $\kappa = \kappa_{U, k}$, we have the following relationship:
\begin{align*}
    \frac{\gamma}{2 - \gamma} &\leq \kappa \leq \gamma , \\
    \implies 0.5\gamma &\leq \kappa \leq \gamma.
\end{align*}
\end{theorem}

\textit{Proof of Theorem \ref{thm:type1totype2}:} we first prove that $\kappa \leq \gamma$. Consider any $L \subseteq U, \abs{S} \leq k$, $L \cap S = \emptyset$. Let $S = \curly{x_1,\dots,x_{\abs{S}}}$. We have
\begin{align*}
    g(L \cup S) - g(L) &= \sum_{j=1}^{\abs{S}} \brac{g\paren{L \cup \curly{x_1,\dots,x_j}} - g\paren{L \cup \curly{x_1,\dots, x_{j-1}}}} \\
    &\leq \sum_{j=1}^{\abs{S}} \frac{1}{\kappa} \brac{g\paren{L\cup \curly{x_j}} - g(L)}, \\
    \implies \kappa &\leq \frac{\sum_{j=1}^{\abs{S}} \brac{g\paren{L\cup \curly{x_j}} - g(L)}}{g(L \cup S) - g(L)}.
\end{align*}
Taking the minimum of the right-hand side of the last equation, we get $\kappa \leq \gamma$ as desired.

We now prove that $\frac{1}{2}\gamma \leq \frac{\gamma}{2 - \gamma} \leq \kappa$. Consider any set $\curly{j,k}$. We have
\begin{align}
g(L \cup \curly{j}) - g(L) &\geq \kappa \paren{g(L \cup \curly{j,k}) - g(L \cup \curly{k})} \notag \\
&= \kappa \paren{g(L \cup \curly{j,k}) - g(L)} -\kappa \paren{g(L \cup \curly{k}) - g(L)} \label{eq:gamma1} \\
g(L \cup \curly{k}) - g(L) &\geq\kappa \paren{g(L \cup \curly{j,k}) - g(L)} -\kappa \paren{g(L \cup \curly{j}) - g(L)} \label{eq:gamma2}
\end{align}

Define
\begin{align*}
    \gamma_1 &= \frac{g(L \cup \curly{j}) - g(L)}{g(L \cup \curly{j,k}) - g(L)} \\
    \gamma_2 &= \frac{g(L \cup \curly{k}) - g(L)}{g(L \cup \curly{j,k}) - g(L)}.
\end{align*}

Rearranging inequalities \ref{eq:gamma1} and \ref{eq:gamma2}, we have
\begin{align*}
    \gamma_1 +\kappa\gamma_2 &\geq\kappa \\
    \gamma_2 +\kappa\gamma_1 &\geq\kappa.
\end{align*}
Furthermore, we have
\[
\gamma_1 + \gamma_2 = \frac{\paren{g(L \cup \curly{j}) - g(L)} + \paren{g(L \cup \curly{k}) - g(L)}}{g(L \cup \curly{j,k}) - g(L)} \geq \gamma.
\]

We consider the following simple minimization problem
\begin{align*}
    \min &\,\,\gamma_1 + \gamma_2 \\
    \text{s.t.} &\,\, \gamma_1 +\kappa\gamma_2 \geq\kappa \\
    &\,\, \gamma_2 +\kappa\gamma_1 \geq\kappa \\
    &\,\, \gamma_1, \gamma_2 \geq 0.
\end{align*}

The optimal value of the above problem is $\frac{2\kappa}{1+\kappa}$, which implies
\[
\gamma \leq \frac{2\kappa}{1+\kappa} \iff\kappa \geq \frac{\gamma}{2 - \gamma}.
\]
One can verify that $0.5\gamma \leq \frac{\gamma}{2 - \gamma}$, which completes the proof of the theorem. $\qedhere$

As we have previously noted, even though there is some intersection between the world of convex sparse atomic optimization and weak submodular maximization, weak submodularity is independently an important notion to ensure the approximation ratios for greedy algorithms. These two alternative conditions lead to similar worst case performance guarantee, however the settings under which they operate may be completely different. To illustrate this point, let us examine below a weakly submodular function which is by no means related to convexity.

Let $h: \: \R^n \to \R$ be a componentwise increasing function, and moreover 
its partial derivatives satisfy $h_i(z) \ge h_i(z') $ for any $z\le z'$ in the domain, where $i=1,2,...,n$.

Such function $h$ clearly exists. For instance, consider a symmetric doubly stochastic matrix $Q$, and define $q(y)=-\frac{1}{2} y^\top Q y +\mathbbm{1}^\top y$ with $\nabla q(y)=-Qy+\mathbbm{1}$ where $\mathbbm{1}$ is the all-ones vector. Let $s(x)=1/(1+\exp(-x))$
be the sigmoid function, with $s'(x)=s(x)(1-s(x))$. Finally, let
\[
h(z)=q(s(z_1),s(z_2),...,s(z_n))
\]
where $z_i\in [0,\infty)$ for $i=1,2,...,n$. Then, for any $0\le z_i\le z'_i$ $(i=1,2,...,n)$ we have
\begin{eqnarray*}
\frac{\partial h}{\partial z_i} (z) &=& \frac{\partial q}{\partial y_i} (s(z)) s(z_i)(1-s(z_i)) \\ 
&\ge & \frac{\partial q}{\partial y_i}(s(z')) s(z'_i)(1-s(z'_i)) = \frac{\partial h}{\partial z_i} (z'),
\end{eqnarray*}
where $i=1,2,...,n$. 
The above inequality holds because if $0\le z\le z'$ then $0<s(z) \le s(z')$, and so $\nabla q(s(z)) - \nabla q(s(z')) \ge 0$; moreover, $s(x)(1-s(x))$ is monotonically decreasing when $x\ge 0$.

We now continue our construction after finding such a function $h$. Let $u_i(t)$ be a unimodal function which attains its maximum at $p_i$ (that is, $u_i(t)$ is increasing for $t<p_i$ and increasing for $t>p_i$). Let us also assume $u_i(0)=0$.

Let $f(x) := h(u_1(x_1), u_1(x_2) ... , u_n(x_n))$. Now we shall show that $g(U)=\max_{x_i,\, i\in U} f(x) - f(0)$ satisfies $(1/c)$-submodularity in the subset sense.

Consider any $U\subset V \subseteq \{1,2,...,n\}$ and $i \not \in V$. Without losing generality, let us denote
\[
U=\{1,2,...,m\}, \, V=\{1,2,...,m,m+1,...,m+\ell\}, \mbox{ and } i=m+\ell+1.
\]
Clearly,
\begin{eqnarray*}
& & g(U\cup \{i\}) - g(U) \\
&=& h(u_1(p_1),...,u_m(p_m),0,...,0,u_{m+\ell+1}(p_{m+\ell+1}),0,...,0) - h(u_1(p_1),...,u_m(p_m),0,...,0) \\
&=& \int_0^{u_{m+\ell+1}(p_{m+\ell+1})} \frac{ \partial h }{ \partial x_{m+\ell+1}}( u_1(p_1),...,u_m(p_m),0,...,0,t,0,...,0)dt.
\end{eqnarray*}
Similarly,
\begin{eqnarray*}
& & g(V\cup \{i\}) - g(V) \\
&=& h(u_1(p_1),...,u_m(p_m),u_{m+1}(p_{m+1}),...,u_{m+\ell}(p_{m+\ell}) ,u_{m+\ell+1}(p_{m+\ell+1}),0,...,0) \\
& & - h(u_1(p_1),...,u_m(p_m),u_{m+1}(p_{m+1}),...,u_{m+\ell}(p_{m+\ell}),0,...,0) \\
&=& \int_0^{u_{m+\ell+1}(p_{m+\ell+1})} \frac{ \partial h }{ \partial x_{m+\ell+1}}( u_1(p_1),...,u_m(p_m),u_{m+1}(p_{m+1}),...,u_{m+\ell}(p_{m+\ell}),t,0,...,0)dt \\
&\le & \int_0^{u_{m+\ell+1}(p_{m+\ell+1})} \frac{ \partial h }{ \partial x_{m+\ell+1}}( u_1(p_1),...,u_m(p_m),0,...,0,t,0,...,0)dt \\
&=& g(U\cup \{i\}) - g(U) .
\end{eqnarray*}

Therefore, the subset submodularity holds with $\kappa=1$. One notices that there is no concavity, restricted or not, required at all in this example.

We now derive the performance of the greedy algorithm with respect to the subset submodularity ratio. Consider the following maximization problem:
\begin{align*}
    \max &\quad g(S) \\
    \text{s.t.} &\quad \abs{S} \leq k, \quad S \subseteq [n],
\end{align*}

where we make some assumptions on $g$.
\begin{enumerate}
    \item Monotonicity: if $U \subseteq V \subseteq [n]$, then $g(U) \leq g(V)$.

    \item Scaled: $g(\emptyset) = 0$.
\end{enumerate}

\begin{algorithm}
\caption{Greedy$(g, r)$}
\begin{algorithmic}[1]
\State $S_0 \leftarrow \emptyset$
\For{$t = 0,\dots,r-1$}
\State $s_{t+1} \leftarrow \argmax_{i \not\in S_t} g\paren{S_t \cup \curly{i}}$
\State $S_{t+1} \leftarrow S_t \cup \curly{s_{t+1}}$
\State $t \leftarrow t + 1$
\EndFor
\State \Return $S_r, g(S_r)$
\end{algorithmic}
\end{algorithm}

\begin{theorem} \label{thm:wksub2}
Let $S^* = \argmax_{\abs{S} = r} g(S)$. At iteration $i$ of Greedy$(g,r)$, we have
\[
g\paren{S_i} \geq \paren{1 - \paren{1 - \frac{\kappa}{r}}^i}g(S^*),
\]
where $\kappa$ is the subset submodularity ratio. In particular, after the greedy algorithm terminates at step $r$, we have
\begin{align*}
    g\paren{S_r} &\geq \paren{1 - \paren{1 - \frac{\kappa}{r}}^r}g(S^*) \\
    &\geq \paren{1 - \exp(-\kappa)} g(S^*).
\end{align*}
\end{theorem}

\textit{Proof of Theorem \ref{thm:wksub2}:} Let $S^* = \curly{x_1,\dots, x_r}$. By monotonicity, we have
\begin{align*}
    g(S^*) &\leq g\paren{S^* \cup S_i} \\
    &= g(S_i) + \paren{g\paren{S_i \cup \curly{x_1}} - g(S_i)} +  \sum_{j=2}^r \paren{g\paren{S_i \cup \curly{x_1,\dots,x_{j-1},x_j}} - g\paren{S_i \cup \curly{x_1,\dots, x_{j-1}}}} \\
    &\leq g(S_i) + \paren{g(S_{i+1}) - g(S_i)} + \sum_{j=2}^r \frac{1}{\kappa}
    \paren{g(S_{i} \cup \curly{x_j}) - g(S_i)} \\
    &\leq g(S_i) + \paren{g(S_{i+1}) - g(S_i)} + \frac{r-1}{\kappa} \paren{g(S_{i+1}) - g(S_i)} \\
    &\leq g(S_i) + \frac{r}{\kappa} \paren{g(S_{i+1}) - g(S_i)}.
\end{align*}

To establish the approximation guarantee, we use induction. The base case $i = 0$ is trivial:
\[
g(S_0) = 0 \geq \paren{1 - \paren{1 - \frac{\kappa}{r}}^0}g(S^*) = 0.
\]
Assume the induction hypothesis
\[
g(S_i) \geq \paren{1 - \paren{1 - \frac{\kappa}{r}}^i}g(S^*).
\]
From the inequality we established earlier, we have
\begin{align*}
    g(S^*) &\leq g(S_i) + \frac{r}{\kappa} \paren{g(S_{i+1}) - g(S_i)} \\
    \iff g(S_{i+1}) &\geq g(S_i) + \frac{\kappa}{r} \paren{g(S^*) - g(S_i)} \\
    &= \frac{\kappa}{r} g(S^*) + \paren{1 - \frac{\kappa}{r}}g(S_i) \\
    &\geq \frac{\kappa}{r} g(S^*) + \paren{1 - \frac{\kappa}{r}}\paren{1 - \paren{1 - \frac{\kappa}{r}}^i} g(S^*) \\
    &\geq \paren{1 - \paren{1 - \frac{\kappa}{r}}^{i+1}}g(S^*).
\end{align*}
This completes the induction. Setting $i = r$, we get the desired approximation guarantee.
$\qedhere$

\section{Discussion and Remarks}

In this paper, we borrowed the notion of atomic sets to capture a general notion of sparsity. We proved that for (restricted) strongly convex, smooth objectives, the performance of the greedy algorithm for unconstrained optimization is connected to a geometric quantity of the objective, the condition number $\mu / L$, and a geometric quantity of the atomic set, $\theta_{r}$. By deriving the $\theta_{r}$ values for various atomic sets, we established explicit approximation guarantees for greedy-type algorithms in various settings. We recovered the ``strong'' approximation guarantees appearing in recent literature for the feature selection and low-rank matrix optimization settings, where ``strong'' is quantified by the fact that the greedy algorithm will find a constant-factor approximation of the optimal $r$-sparse solution within $O(r)$ iterations. Namely, we can guarantee the greedy algorithm will produce a \textit{sparse} solution that attains objective value on the same order as the optimal solution. Through the $\theta_{r}$ value, we have provided a simple method of deriving greedy approximation guarantees for theoretically any atomic set. We believe that the $\theta_{r}$ value of an atomic set furthermore holds some degree of truth regarding the ability of greedy search to discover good sparse solutions, since it measures the ability of the atomic set to approximate any element of an arbitrary $r$-subspace. Finding the atom best aligned to the gradient of the objective function is precisely the ``greedy'' part of one variant of the greedy algorithm, which have previously mentioned to be the computationally feasible alternative of pure greedy selection.

A secondary contribution of this paper was disentangling the notion of weak submodularity from the performance of greedy algorithms for continuous, convex optimization. Certainly, strongly convex, smooth objective functions on certain atomic sets can be shown to satisfy explicit degrees of weak submodularity. However, two key properties required of an atomic set to fruitfully connect weak submodularity and convex optimization is first whether the atomic set contains (approximately) orthogonal bases to the ambient space, and second whether the union of two sparse sets of atoms can be sparsely reparameterized into a set of (approximately) orthogonal atoms such that the span of the former is equal to (or contained in) the span of the latter. If the atomic set cannot guarantee the above two properties, one is hard-pressed to leverage the properties of strong convexity and smoothness to bound weak submodularity. However, by separately proving approximation guarantees of greedy algorithms for strongly convex, smooth functions versus weakly submodular functions, we now have access to guarantees of greedy algorithms on a much richer variety of functions.

Since our model is unconstrained optimization, we note that our model does not account for certain atomic sets where the norm of the atoms comes into play. For example, take the problem of recovering a low-rank and sparse matrix decomposition, such as in robust PCA \cite{CLMW11}. Earlier literature has shown that this problem admits a convex relaxation. Namely, one can frame it as an atomic norm regularization problem, where the corresponding atomic norm ball is the convex hull of rank-one matrices $\curly{uv^\top}$ and some multiple of the singleton matrices $\gamma \curly{e_ie_j^\top}$ \cite{CSPW11}. However, our greedy algorithm is agnostic to the norm of atoms, and since singleton matrices are themselves rank-one, optimizing over the above atomic set is equivalent to optimizing over the normalized rank-one matrices. In general, this issue occurs when one wants to promote selecting from a subfamily of a larger atomic set. This is fundamental to the unconstrained nature of the optimization, and thus there is no obvious adjustment that can be made.

An possibly interesting direction of future research would be to address the sparsity parameter in the greedy algorithm. As stated in this paper, the sparsity parameter is entirely user-determined, and thus it is up to the user to somehow determine the ``true'' sparsity of the model. It may be interesting to consider principled, yet computationally cheap approaches to selecting an ``optimal'' sparsity. As an example, we previously mentioned that PCA is a special case of applying the greedy algorithm to the linear recovery objective. In that case, parallel analysis \cite{Horn65, BE92} has long been known as a very good heuristic for estimating the rank, and much more recently been shown theoretically to select the optimal rank under some model and noise assumptions \cite{Dob17, DO17}. Heuristics and theory either for more general objective functions or different atomic sets will likely involve vastly different tools, but it is nevertheless an interesting problem.

\textbf{Acknowledgements:} The author is currently a postgraduate researcher in the lab of Prof.\ Yuval Kluger, Yale Applied Mathematics Program, and would like to thank Prof. Kluger for his support. The author would also like to thank Prof.\ Sahand Negahban for introducing him to the topic as an undergraduate thesis advisor, as well as many helpful comments on a nascent version of this paper. Any errors and aberrant conclusions are the author's own.

\bibliography{refs}{}
\bibliographystyle{ieeetr}

\appendix

\section{Computational Complexity of the Atomic Condition Number}\label{apdx:condnumber}

An immediate question one might ask is: given an arbitrary atomic set $\A$, can we numerically compute $\theta(\A)$ efficiently? As one may suspect, the atomic condition number cannot be computed efficiently in general. It suffices to restrict our attention to $\A \subseteq \R^n$ and $\abs{\A} = m < \infty$. For convenience, let us assume $\norm{a_i}_2 = 1$ for all $a_i \in \A$. Let $A \in \R^{n \times m}$ denote the matrix whose columns are $a_i \in \A$. Consider the following definition of $\theta$:
\begin{align*}
    \theta &= \min_{\norm{x}_2 = 1} \max_{a \in \A} \abs{\ip{x, a}} \\
    &= \min_{x \in \R^n} \frac{\norm{A^\top x}_\infty}{\norm{x}_2}.
\end{align*}

Observe that when $A$ is degenerate, i.e.\ $\A$ does not span $\R^n$, then $\theta = 0$. Therefore, we need only to consider $A$ such that $m \geq n$ and $\rank{A} \geq n$. Observe that if we can solve the following problem efficiently, we can also compute $\theta$ efficiently:
\begin{align*}
    (P)\quad \max\,\,& \norm{x}_2^2 \\
    \text{s.t.}\,\,& A^\top x \leq t\ones \\
    &A^\top x \geq -t\ones.
\end{align*}
Setting $y = A^\top x$ and $Q \defeq (A A^\top)^{-1} \succ 0$, we may re-write $(P)$ as
\begin{align*}
    (P')\quad \max\,\,& y^\top Q y \\
    \text{s.t.}\,\,& -1 \leq y_i \leq 1 \quad \forall i \in [m].
\end{align*}
Since $Q$ is positive-definite, $(P')$ precisely coincides with a sub-case of the \textsc{MaxQP} problem considered in \cite{CW04}, where it is shown that quadratic programming over the $\ell^\infty$ hypercube is NP-hard, as the \textsc{MaxCut} problem is a special case \textsc{MaxQP}. Therefore, precisely computing $\theta$ is difficult in general. However, this does not preclude the possibility of estimating $\theta$ efficiently either through approximation schemes or local-search heuristics: estimation is certainly possible when $\A \subset \R^n$ and $\abs{\A} < \infty$, but may not be as straightforward when $\A$ is uncountably infinite or a large combinatorial set, such as the cases of the rank-one matrices or the permutation matrices.

\section{Illustrative Experiments}

We supplement our theoretical analyses with some numerical experiments. We note that there are plentiful experiments involving greedy algorithms applied on various large-scale real data for the more popular atomic sets, for example unit basis vectors \cite{KEDN17sparse, EKDN18, DK11} and low-rank matrices \cite{KEDN17, SGS11}. In the first experiment, we measure the greedy algorithm's performance on linear recovery tasks. That is, we generate a ground truth vector that is a linear combination of a predetermined (sparse) number of atoms, where the weights are $\pm 1$ Bernoulli random variables. We feed to the greedy algorithm the least-squares objective function, where the target vector has been corrupted by ~10\% i.i.d.\ additive gaussian noise, as well as an overestimate of the sparsity. Since the realization of each iterate of the greedy algorithm does not depend on the sparsity parameter, one can easily determine the performance of the greedy algorithm had the sparsity been under, perfectly, or over-estimated. We validate the performance of the greedy algorithm on a different noisy measurement of the target vector. The lines indicate the true sparsity of the solution, as well as the objective/validation functions evaluated on the true target vector.

\begin{figure}[h]
\centering
\includegraphics[scale = 0.42]{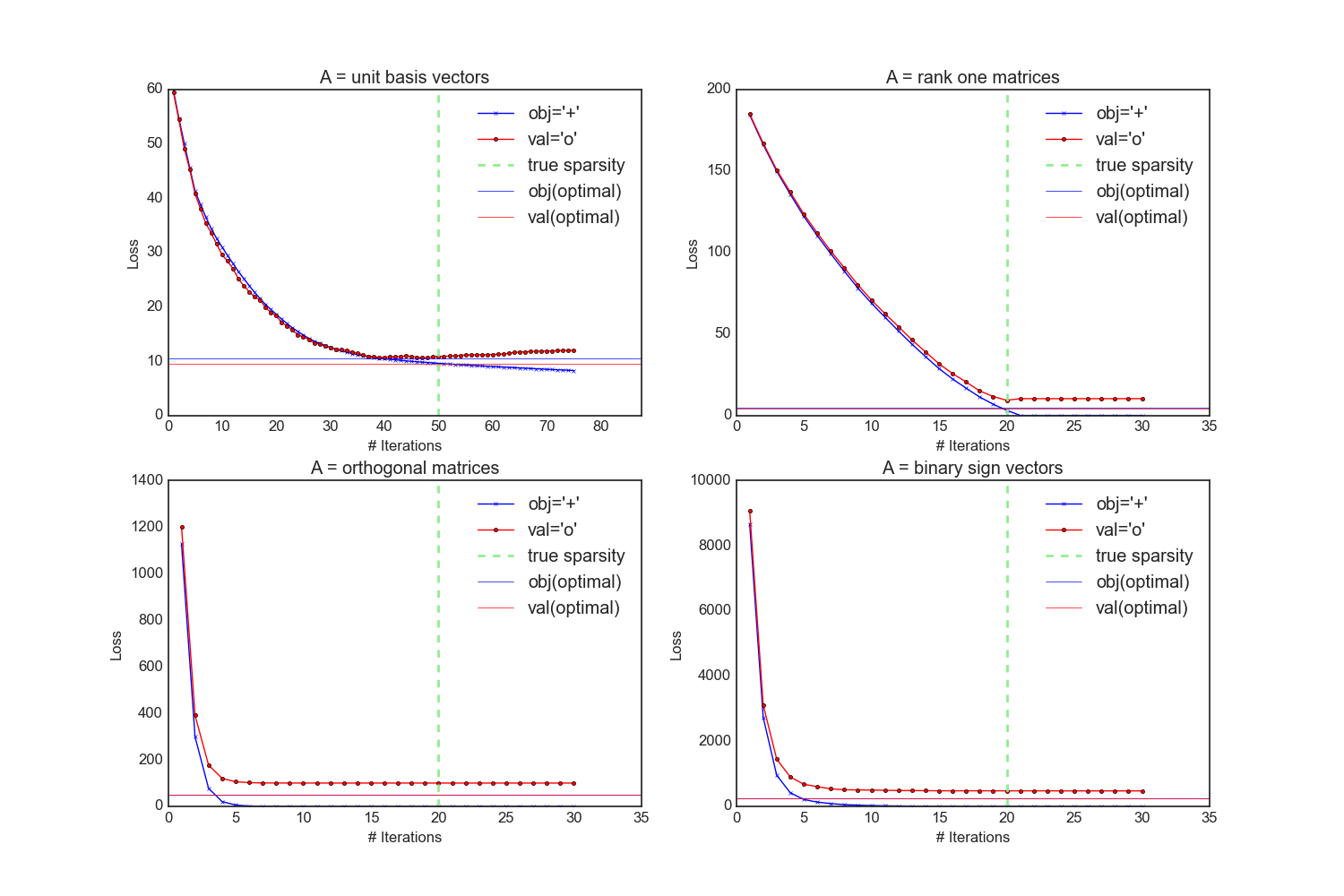}
\caption{Performance of the greedy algorithm on linear recovery tasks for various atomic sets}
\end{figure}

We observe that the greedy algorithm performs extremely well for all the above recovery tasks, and finds optimal sparse solutions for all the above atomic sets. This is unsurprising for the unit basis vectors and rank-one matrices: for the unit basis vectors, each iteration of the greedy algorithm finds the next largest entry in absolute value; for the rank-one matrices, the greedy algorithm is equivalent to rank-$k$ PCA or SVD of the measurement matrix, where $k$ is the inputted sparsity parameter. Both these schemes will essentially find the optimal solution under our noise assumptions. What is slightly surprising is that the greedy algorithm converges to optimality in the orthogonal matrix and sign vectors case, especially how quickly it converges to $0$. One possible explanation of this convergence rate might be related to the identifiability of linear combinations of those atoms. While in the unit basis and rank-one matrices case, the linear combination of $k$ distinct atoms will (almost always) at best be $k$-sparse, in the orthogonal matrix and sign vectors case identifiability issues may arise, where certain linear combinations of $k$ atoms are actually better represented by fewer atoms.

We demonstrate the relationship between certain quantities of an atomic set and the atomic condition number. Given an arbitrary set of vectors as an atomic set, it is a difficult problem to estimate the atomic condition number, as demonstrated in section \ref{apdx:condnumber} of the appendix. We restrict our attention to collection of unit vectors in $\R^n$, where computing the atomic condition number reduces to solving
\begin{align*}
    \theta &= \min_{\norm{x}_2 = 1} \max_{a \in \A} \abs{\ip{x, a}} \\
    &= \min_{x \in \R^n} \frac{\norm{A^\top x}_\infty}{\norm{x}_2},
\end{align*}
where $A$ is the matrix that contains the atoms as column vectors. Given $A$, the latter value can be efficiently estimated in practice using local-search methods. As a preliminary, we observe that $\theta = n^{-1/2}$ is maximal for $A \in M_n$, and is attained if and only if $A$ is orthogonal. This can be easily established by observing that, given an arbitrary non-singular matrix $A$, $\theta$ is monotonically increased by taking a vector from $A$ and orthogonalizing it against the rest of the atoms. In the first experiment, we start with $\A = \curly{e_i}_{i=1}^5$, $A = I_5$. We then corrupt the first column of $A$ by setting $e_1 = v/\norm{v}$, where $v = \lambda e_1 + (1 - \lambda)s$ for different values of $0 \leq \lambda \leq 1$, where $s \in \curly{\pm 1}^5$ is a fixed binary sign vector and $s(1) = -1$. For the latter two experiments, we generate many $A$ randomly as i.i.d.\ gaussian random matrices and normalize the columns, then for each $A$ we measure two values that are somewhat related to how close to ''orthogonal'' the matrix is, namely the mean coherence (mean of $\abs{\ip{a_i, a_j}}$ for each $i \neq j$) and smallest singular value of $A$ and plot the relationship with the atomic condition number.

\begin{figure}[h]
\centering
\includegraphics[scale = 0.42]{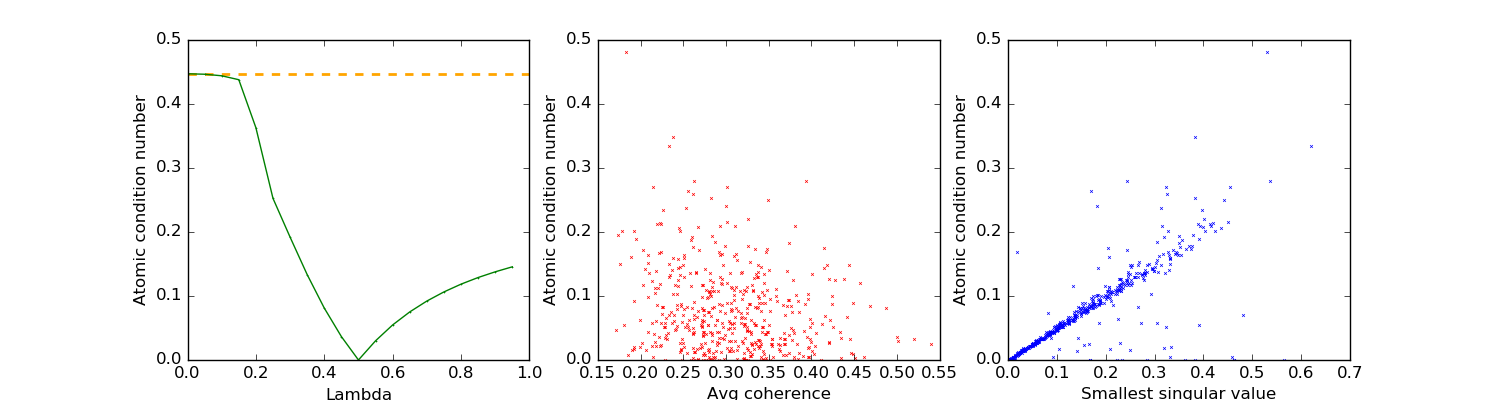}
\caption{Relationship of the atomic condition number $\theta$ and various values}
\end{figure}

In the leftmost graph corresponding to the first experiment, the dotted orange line indicates the optimal $\theta = 1/\sqrt{5}$. The trend observed is not particularly surprising. At $\lambda = 0.5$, observe that since $s(1) = -1$, the first element of $v = e_1 + s$ is $0$, which means that $A$ is singular, and therefore the atomic condition number is $0$. In the plot of atomic condition number versus the mean coherence, we observe a somewhat negative correlation, which makes sense, as a high average coherence indicates that many of the vectors in $A$ are highly correlated, which implies $\theta$ is not optimal. Perhaps most surprising is the plot of atomic condition number versus the smallest singular value of $A$, where one observes an extremely tight linear relationship between the smallest singular value and atomic condition number (where many of the outliers can be attributed to the numerical inconsistency of the local search heuristic used to compute $\theta$). At first glance, this might seem feasible, as the smallest singular value and $\theta$ can be written in similar ways:
\begin{align*}
    \sigma_{\min} = \min_{x \in \R^n} \frac{\norm{Ax}_2}{\norm{x}_2}, \quad \theta = \min_{x \in \R^n} \frac{\norm{Ax}_\infty}{\norm{x}_2}.
\end{align*}
However, there cannot be a linear relationship between the smallest singular value of $A$ and the atomic condition number, as we have previously established in Appendix \ref{apdx:condnumber} that computing the atomic condition number of $\A \subset \R^n$ is NP-hard, while computing the smallest singular value is not.

\section{Computing Sparse Atomic Condition Numbers in Table~\ref{sparse theta values}}


We recall the mathematical definitions of the atomic and sparse atomic condition numbers $\theta$ and $\theta_r$: given the ambient vector space $V$ of atomic set $\A$
\begin{align*}
    \theta &\defeq \min_{\norm{v} = 1} \max_{a \in \A} \frac{\abs{\ip{v, a}}}{\norm{a}} \\
    \theta_r &\defeq \min_{\substack{U \subset \A \\ \abs{U} \leq r}} \min_{v \in \Span{U}} \max_{a \in \A} \frac{\abs{\ip{v, a}}}{\norm{a}}.
\end{align*}
$\theta$ measures the ability of the atomic set to approximate any vector in the ambient space, where $\theta_r$ measures the ability to approximate any vector given it comes from a subspace of dimension no more than $r$. We observe that $\theta_r \geq \theta$, and $\theta_r \geq \theta_k$ if $r \geq k$.


\subsubsection*{Standard Basis}
The derivations of $\theta$ and $\theta_r$ are essentially identical. For $\theta$, given any $v \in \R^n$, $\max_{i \in [n]} \abs{\ip{v, e_i}}/\norm{v}$ is attained by the largest entry in absolute value. It is then easy to see that
\[
\theta = \min_{\norm{v} = 1} \max_{i \in [n]} \abs{\ip{v, e_i}}
\]
is attained by the scaled all-ones vector: $v = n^{-1/2} \ones$. Similarly, if $U \subset \curly{e_i}^n$, $\abs{U} = r$, then $\theta_r$ is attained by the scaled indicator vector of $U$: $v = r^{-1/2} \ones_U$, which gives us $\theta_r = r^{-1/2}$. We note that by finding the vectors that attain $\theta$ and $\theta_r$, our lower bounds are tight.

The $\theta$ and $\theta_r$ values for the standard basis hold for any orthogonal basis, a fact that we will use in the analysis of the 2-ortho basis case.

\subsubsection*{Rank-One Matrices}
The derivation of $\theta$ and $\theta_r$ for $\A = \curly{uv^\top}$ is very similar to that for the standard basis, where we instead look at the spectrum of a given matrix. For $\theta$, if $M \in \R^{m,n}$, we have $\rank{M} \leq \min\curly{m,n}$.
\begin{align*}
    \theta &= \min_M \max_{uv^\top, \norm{u}\norm{v}=1} \frac{\ip{M, uv^\top}}{\norm{M}_F} \\
    &= \min_M \max_{uv^\top, \norm{u}\norm{v}=1} \frac{1}{\norm{M}_F} \tr{u^\top M v} \\
    &= \min_M \frac{1}{\norm{M}_F} \sigma_1(M)
\end{align*}
where $\sigma_1(M)$ is the largest singular value of $M$. Defining $\sigma(M)$ as the vector containing the singular values of $M$, we recall that $\norm{M}_F = \norm{\sigma(M)}_2$. Therefore, as in the case of the standard basis, the spectrum vector that attains $\theta$ is the scaled all-ones: $\sigma(M) = n^{-1/2} \ones$, which of course has leading singular value $\sigma_1(M) = n^{-1/2}$.

For $\theta_r$, given $U \subset \curly{uv^\top}$, $\abs{U} = r$, any matrix $M \in \Span{V}$ is the linear combination of at most $r$ rank-one matrices, and therefore is at most rank $r$. In other words, we have $\abs{\supp{\sigma(M)}} \leq r$. Following the argument for the standard basis, we have $\theta = r^{-1/2}$. The values derived for $\theta$ and $\theta_r$ are tight.

\subsubsection*{Disjoint Group Sparse Atoms}

Given the standard basis $\curly{e_i}^n \subset \R^n$, disjoint group-sparse atoms are defined as the elements of a partition of the basis: $\A = \mathcal{P}(\curly{e_i}^n)$. Let $\abs{\A} = L$. Note that these atoms are not vectors, and thus we must make a few adjustments to some definitions. Given objective function $f$, the accompanying set function $g$ is now defined
\[
g(\mathcal{S}) = \max_{x \in \Span{\bigcup_{P \in \mathcal{S}} P}} f(x) - f(0),
\]
where $\mathcal{S} \subset \mathcal{P}$ and $P$ are the groups contained in $\mathcal{S}$. In other words, $g(\mathcal{S})$ returns the shifted optimal value of $f$ searching over the vectors whose support lies in the union of the groups. Applying this new definition to our algorithm, $\texttt{PureGreedy}(\A, f, L_{t-1}, \beta)$ makes sense as is. For $\texttt{OMPSel}(\A, f, L_{t-1}, \beta)$, we modify the definition such that it returns the group $P_t$ satisfying
\[
\norm{\proj_{P_t}\paren{\nabla f\paren{B^{(L_{t-1})}}}} \geq \beta \max_{P \in \A} \norm{\proj_{P}\paren{\nabla f\paren{B^{(L_{t-1})}}}}.
\]
In other words, whereas usually $\texttt{OMPSel}$ finds the atom (vector) that best explains the gradient of the previous iterate, it now finds the group (subspace) that can best explain the gradient of the previous iterate. We are now able to define $\theta$ and $\theta_r$.
\begin{align*}
    \theta &= \min_{\norm{v} = 1} \max_{P \in \A} \frac{\abs{\ip{v, \proj_P(v)}}}{\norm{\proj_P(v)}} \\
    &= \min_{\norm{v} = 1} \max_{P \in \A} \norm{\proj_P(v)},
\end{align*}
and similarly,
\begin{align*}
    \theta_r &= \min_{\substack{U \subset \A \\ \abs{U} = r}} \min_{\substack{v \in \Span{\bigcup_{P \in U} P} \\ \norm{v} = 1}} \max_{P \in \A} \norm{\proj_P(v)}.
\end{align*}

It suffices for us to lower bound $\theta_r$; the proof for $\theta$ is identical, setting $r = L$. Given any partition $\mathcal{P}$ of $\R^n$, and any subset $U \subseteq \mathcal{P}$, $\abs{U} = r$, we consider what the minimizing vector $v$ that attains $\theta_U$ would look like. Analogous to what we saw in the proofs for the standard basis and rank-one matrices, if we can find a vector $v$ such that $\norm{\proj_P(v)} = r^{-1/2}$ for all $P \in U$, then $v$ is clearly the minimizing vector that attains $\theta_U$, since $\norm{v}=1$ and shifting any mass around in the vector can only cause $\max_{P \in U} \norm{\proj_P(v)}$ to increase. Let us construct such a $v$. Consider a vector $v = \sum_{i=1}^r c_i \ones_{P_i}$, where $\ones_{P_i} \in \R^n$ is the indicator vector corresponding to the group $P_i \in U$. We have the constraint $\norm{v} = 1$, and the property $\norm{\proj_{P_i}(v)} = r^{-1/2}$ for all $i$. In other words, a vector $v$ attaining $\theta_U$ would satisfy
\begin{align*}
    \sum_{i=1}^r c_i^2 \abs{P_i} &= 1 \\
    \abs{c_i}\sqrt{\abs{P_i}} &= r^{-1/2} \quad i = 1,\dots,r.
\end{align*}
A possible solution is $c_i = \paren{r\abs{P_i}}^{-1/2},\,\,i=1,\dots,r$. Substituting $v$ into our expression for $\theta$ and $\theta_r$, we get $\theta = n^{-1/2}$ and $\theta_r = r^{-1/2}$, and these lower bounds are tight our construction.

\subsubsection*{2-Ortho Basis}

In the signal processing and dictionary learning communities, the 2-ortho basis refers to the union of the standard and Fourier orthogonal bases \cite{CDS01, Elad10}. We will consider a general union of two orthogonal bases, w.l.o.g.\ the union of the standard basis and an arbitrary orthogonal basis. Let us denote this atomic set
\[
\A \defeq \curly{e_i}_{i=1}^n \cup \curly{\psi_j}_{j=1}^n \subset \R^n.
\]
A popular, albeit crude, measure of distance between bases is the ``mutual coherence''. For example, we have for the standard-union-Fourier,
\[
\mu(\A) \defeq \max_{i,j} \abs{e_i^\top \psi_j} = n^{-1/2}.
\]
We note that $n^{-1/2}$ is the lower bound on the mutual coherence of two orthogonal bases. We also recall that the (non-sparse) atomic condition number $\theta$ is monotonically increasing with respect to adding more elements. One might therefore hope that there is some link between the mutual coherence and a factor of ``improvement'' to $\theta$. However, we will show that in general, the atomic condition number $\theta$ cannot be improved above $n^{-1/2}$ by example.

We consider the extremal case of the Hadamard basis (matrix), which is an orthogonal basis consisting of vectors in $\curly{\pm 1}^n$. To be sure, Hadamard bases do not exist for every dimension (consider any odd dimension greater than $1$), but there exist infinitely many. We consider a special infinite subfamily of the Hadamard bases: the regular Hadamard bases. Regular Hadamard bases are simply Hadamard matrices whose row and column sums are all equal, which restricts the dimension of the basis to square numbers. In particular, if the dimension of the regular Hadamard basis/matrix is $n = 4u^2$, then the row and column sums are all equal to $\pm 2u$, which further implies that each column has $2u^2 \pm u$ positive entries and $2u^2 \mp u$ negative entries. It is simple to see that any Hadamard basis attains the minimal $n^{-1/2}$ mutual coherence with the standard basis. We now consider the scaled all-ones vector $n^{-1/2}\ones$, which we recall attains $\theta$ when $\A$ is the standard basis. For any member $v$ of the regular Hadamard basis, its inner product with the all-ones vector is also
\[
\frac{\abs{\ip{v, \ones}}}{\norm{v}\norm{\ones}} = \frac{1}{n}2u = \frac{1}{n}\sqrt{n} = n^{-1/2}.
\]
In other words, the all-ones vector \textit{also} attains $\theta$ for the regular Hadamard basis. Therefore, $\theta$ of the union of the standard and Hadamard basis is still $n^{-1/2}$. This means that in full generality, there may be no relation between the mutual coherence of two orthogonal bases $\mu\paren{\A}$ and the atomic condition number $\theta$.

It remains to be seen what happens to the sparse atomic condition number $\theta_r$. We have seen that arbitrarily adding atoms to the standard basis, though improving $\theta$, may instantaneously cause $\theta_r$ to shrink from $r^{-1/2}$ to $n^{-1/2}$. Therefore, the more may not be the merrier when it comes to the sparse atomic condition number. However, we note that the example we used to demonstrate corrupting the $\theta_r$ value, where we added $e_1 + \varepsilon \ones$ to $\curly{e_i}_{i=1}^n$, would have a mutual coherence arbitrarily close to $1$ (since we are essentially adding a slightly perturbed member of the standard basis). However, what happens when we can guarantee two orthogonal bases are a certain angle away from each other? We will show that when the sparsity is below a certain level $r \leq \mu\paren{\A}^{-1}$, then the sparse atomic condition number is lower bounded by $\Omega\paren{r^{-1/2}}$. In other words, under sufficient sparsity, the sparse condition number of a 2-ortho basis is on the same order as just one orthogonal basis.

Let $\curly{e_i}_{i=1}^n$ be an orthogonal basis (w.l.o.g.\ the standard basis) and $\curly{\psi_j}_{j=1}^n$ be another orthogonal basis, with mutual coherence $\mu$. Consider any vector
\[
x = \sum_{i \in I} c_i e_i + \sum_{j \in J} d_j \psi_j \neq 0
\]
where $\abs{I} = r_1$, $\abs{J} = r_2$, $r_1 + r_2 \leq \mu^{-1}$. For our purposes, we can assume $\norm{x} = 1$. We may also assume $r_1, r_2 > 0$, since if either one is $0$, then we are reduced to computing the sparse atomic number of a single orthogonal basis. From the AM-GM inequality, we have $\sqrt{r_1 r_2} \leq (r_1 + r_2)/2 \leq \mu^{-1}/2$. We also have
\begin{align*}
    \norm{x}^2 &= \norm{c}^2 + \sum_{i \in I} \sum_{j \in J} c_i d_j \ip{e_i, \psi_j} + \norm{d}^2,
\end{align*}
where $c$ and $d$ are the coefficient vectors from $x$. From the definition of mutual coherence, we have $-\mu^{-1} \leq \ip{e_i, \psi_j} \leq \mu^{-1}$. Thus, we can make a crude upper and lower bound on $\norm{x}^2$ with respect to $\norm{c}$ and $\norm{d}$:
\begin{align*}
    \norm{x}^2 &\leq \norm{c}^2 + \norm{d}^2 + \norm{c}_1\norm{d}_1 \mu \\
    &\leq \norm{c}^2 + \norm{d}^2 + \sqrt{r_1 r_2} \norm{c}\norm{d}\mu \\
    &\leq \norm{c}^2 + \norm{d}^2 + \frac{1}{2}\norm{c}\norm{d} \\
    &\leq \frac{5}{4}\paren{\norm{c}^2 + \norm{d}^2} , \\
    \norm{x}^2 &\geq \norm{c}^2 + \norm{d}^2 - \norm{c}_1\norm{d}_1 \mu \\
    &\geq \frac{3}{4}\paren{\norm{c}^2 + \norm{d}^2}.
\end{align*}
Without loss of generality, let $\norm{c} \geq \norm{d}$. Hence, $\norm{c} \geq \sqrt{\frac{4}{10}} \norm{x}$. Let $c_p = \argmax_{i \in I} \abs{c_i}$. Letting $\norm{x}= 1$, we have
\begin{align*}
    \frac{\abs{\ip{x, c_p}}}{\norm{x} \norm{c_p}} &\geq \abs{c_p} - \mu \sum_{j \in J} \abs{d_j} \\
    &\geq \norm{c}/\sqrt{r_1} - \sqrt{r_2} \mu \norm{d} \\
    &\geq \norm{c}/\sqrt{r_1} - \norm{d}/\paren{2\sqrt{r_1}} \quad \mbox{(since $\sqrt{r_1r_2} \leq \mu^{-1}/2$)} \\
    &\geq \norm{c}/\paren{2\sqrt{r_1}} \\
    &\geq \norm{c}/\paren{2\sqrt{r}} \\
    &\geq \paren{\sqrt{4/10} \norm{x}} \frac{1}{2} r^{-1/2} \\
    &\geq 10^{-1/2} r^{-1/2} \\
    &= \Omega\paren{r^{-1/2}} . 
\end{align*}
Therefore, for any $r \leq \mu^{-1}$, we have
\begin{align*}
    \theta_r &\defeq \min_{\substack{U \subset \A \\ \abs{U} \leq r}} \min_{v \in \Span{U}} \max_{a \in \A} \frac{\abs{\ip{v, a}}}{\norm{a}} \\
    &= \Omega\paren{r^{-1/2}}.
\end{align*}
We note that $\mu^{-1}$ is at maximum $\sqrt{n}$, which is attained by the classic standard-union-Fourier basis or the standard-union-Hadamard basis we discussed.

\subsubsection*{Binary sign vectors}

Given vector $v$, it is clear that the atom $a \in \curly{\pm 1}^n$ that maximizes $\frac{\ip{v, a}}{\norm{a}}$ is $a = \text{sign}(v)$, where $\text{sign}(v)$ is the binary sign vector whose entries are the signs of the entries of $v$. Therefore, we have
\begin{align*}
    \theta &= \min_{\norm{v} = 1} \max_{a \in \A} \frac{\abs{\ip{v,a}}}{\norm{a}} \\
    &= \min_{\norm{v} = 1} n^{-1/2} \ip{v, \text{ sign}(v)} \\
    &= \min_{\norm{v} = 1} n^{-1/2} \norm{v}_1 \\
    &= n^{-1/2},
\end{align*}
where the last line comes from the $\ell^1$-$\ell^2$ equivalence of norms inequality, which is attained by $v = e_i$ for any standard basis vector $e_i$. We observe that $\theta_r$ cannot be any larger than $n^{-1/2}$ by a simple example. Let $U$ be any set containing $a_1, a_2$, which are two sign vectors that are identical apart from one entry where the sign is flipped, without loss of generality the first entry. Take $x = \frac{1}{2}a_1 + \frac{1}{2} a_2 = e_1$, which is in $\Span{U}$. As we saw earlier, $e_1$ is a vector that attains $\theta$, and therefore $\theta_r = \theta = n^{-1/2}$.

\subsubsection*{Orthogonal Matrices}

The orthogonal matrices are an interesting atomic set where the $\theta$ and $\theta_r$ values are precisely identical. Let $\norm{\cdot}_F$ denote the Frobenius norm of a matrix, and $\norm{\cdot}_\ast$ denote the nuclear norm (or trace norm). Given any matrix $M$, $\norm{M}_F = 1$ and its singular value decomposition $M = U \Sigma V^\top$, we know that the closest orthogonal matrix in Frobenius norm to $M$ is the matrix $Q = UV^\top$ \cite{HJ12}. From the identity
\begin{align*}
\norm{M - Q}_F^2 &= \norm{M}_F^2 + \norm{Q}_F^2 - 2\ip{M,Q} \\
&= \norm{M}_F^2 + n - 2\ip{M, Q}
\end{align*}
we observe the optimal solutions of the following two problems are identical
\[
\min_{Q:\,Q^\top Q = I_n} \norm{M - Q}_F^2 \equiv \max_{Q:\,Q^\top Q = I_n} \ip{M,Q}.
\]
We can now establish the lower bound for $\theta$
\begin{align*}
    \theta &\defeq \min_{\norm{M}_F = 1} \max_{Q \in \A} \frac{\abs{\ip{M, Q}}}{\norm{Q}} \\
    &= \min_{\norm{M}_F = 1} n^{-1/2} \tr{(UV^\top)^\top M} \\
    &= \min_{\norm{M}_F = 1} n^{-1/2} \tr{\Sigma} \\
    &= \min_{\norm{M}_F = 1} n^{-1/2} \norm{M}_{\ast}.
\end{align*}
Since the nuclear norm and Frobenius norms are special cases of Schatten $p$-norms, with $p = 1,2$ respectively, which are defined as $p$-norms on the singular values of a matrix, we have from the equivalence of norms: $\norm{M}_{\ast} \geq \norm{M}_F$, which gets us the lower bound on $\theta$:
\begin{align*}
    \theta &= \min_{\norm{M}_F = 1} n^{-1/2} \norm{M}_{\ast} \\
    &\geq \min_{\norm{M}_F = 1} n^{-1/2} \norm{M}_F \\
    &= n^{-1/2}.
\end{align*}

It is simple to see that $\theta_r$ cannot be better than $n^{-1/2}$ by considering the following example. Let $r = 2$ and $U = \curly{I_n, Q}$, where
\[
Q = \begin{bmatrix} \paren{(d-1)/d}^{-1/2} & -d^{-1/2} &  \\ d^{-1/2} & \paren{(d-1)/d}^{-1/2} & \\ & & I_{n-2} \end{bmatrix}.
\]
and we set $M = 2\sqrt{1 - \sqrt{\frac{d-1}{d}}}\paren{I_n - Q}$. Letting $d$ be arbitrarily large, we see that $\theta_r$ can be made arbitrarily close to $\theta$ for any $r \geq 2$.

\end{document}